\documentclass [12pt]{article}
\usepackage{latexsym,amsthm,amsfonts,amsmath}

\newtheorem{theo}{Theorem}

\newtheorem{prob}{Problem}

\newcommand{\p}{{\partial}}
\newcommand\N{\ensuremath{\mathbb{N}}}

\newcommand\Z{\ensuremath{\mathbb{Z}}}

\def\br{\bigr}
\def\bl{\bigl}
\def\diam{\mathop{\mathrm{diam}}}

\title{Pinched exponential volume growth implies \\  an infinite dimensional isoperimetric inequality}
\author{Itai Benjamini \and Oded Schramm}

\date{}

\begin{document}
\maketitle
\begin{abstract}
Let $G$ be a graph which satisfies $c^{-1}\,a^r\le |B(v,r)|\le c\,a^r$,
for some constants $c,a>1$, every vertex $v$ and every radius $r$.
We prove that this implies the isoperimetric inequality
$|\partial A| \ge C |A| / \log(2+ |A|)$ for some constant $C=C(a,c)$ and every
finite set of vertices $A$.
\end{abstract}

A graph $G=\bigl(V(G),E(G)\bigr)$ has pinched growth $f(r)$ if
there are two constants $0 < c < C < \infty$ so that every ball
$B(v,r)$ of radius $r$ centered around a vertex $v\in V(G)$ satisfies
$$
c \,f(r) < \bl|B(v,r)\br| < C\, f(r)\,.
$$
For example, Cayley graphs and vertex-transitive graphs have pinched growth.

It is easy to come up with an example of a tree for which every
ball satisfies $\bl|B(v,r)\br| \ge  2^{r/2}$, yet there are
arbitrarily large finite subsets of $G$ with one boundary vertex.
For example, start with $\N$ (an infinite one-sided path) and
connect every vertex $n$ to the root of a
 binary tree of depth $n$.  This tree does not have pinched growth.

We will see that, perhaps surprisingly, the additional assumption
of pinched exponential growth
(that is, pinched growth $a^r$, for some $a>1$)
implies an infinite dimensional
isoperimetric inequality. For a set $A\subset V(G)$ of vertices denote
by $\partial A $ the (vertex) boundary of $A$, consisting of
vertices outside of $A$ which have a neighbor in $A$.

\begin{theo}\label{main}
Let $G$ be an infinite graph with pinched growth $a^r$,
where $a>1$.
Then there is a constant $c>0$ such that for every finite set of vertices
$A\subset V(G)$,
\begin{equation}\label{e.ii}
\bl|\partial A\br| \ge c\, |A| / \log(2+ |A|)\,.
\end{equation}
\end{theo}

We say that $G$ satisfies an $s$-dimensional isoperimetric inequality
if there is a $c>0$ such that $\bl|\p A\br|\ge  c\,|A|^{(s-1)/s}$ holds for every
finite $A\subset V(G)$.
Thus,~(\ref{e.ii}) may be considered an
infinite dimensional isoperimetric inequality.

Coulhon and Saloff-Coste \cite{CS} proved that when $G$ is a Cayley graph of an
infinite, finitely-generated, group, the isoperimetric inequality
\begin{equation}\label{e.cs}
\bl|\p A\br| \ge \frac{|A|}{4\,m\,\phi\bl(2|A|\br)}
\end{equation}
holds for every finite $A\subset V(G)$, where $m$ is the number of neighbors
every vertex has and $\phi(n)=\inf\bl\{r\ge 1:|B(v,r)|\ge n\br\}$
(here, $v\in V(G)$ is arbitrary).
This result implies Theorem~\ref{main} for the case where $G$ is a Cayley graph,
even when the upper bound in the pinched growth condition is dropped.
The tree example discussed above shows that Theorem~\ref{main} is not
valid without the upper bound.  Thus, the (short and elegant)
proof of~(\ref{e.cs}) from \cite{CS} does not generalize to give
Theorem~\ref{main}, and, in fact, the proof below does not seem related to
the arguments from~\cite{CS}.

It is worthwhile to note that~(\ref{e.cs}) is also interesting for
Cayley graphs with sub-exponential growth.  For example, it shows
that $\Z^d$ satisfies a $d$-dimensional isoperimetric inequality.

Another related result, with some remote similarity in the proof,
is due to Babai and Szegedy \cite{BaS}. They prove that for a
finite vertex transitive graph $G$, and $A \subset V(G)$, $0 < |A| <
|G|/2$,
\begin{equation*}
  |\p A\br|  \ge {|A|} /(1+\diam G)\,.
\end{equation*}

\medskip

The isoperimetric inequality~(\ref{e.ii}) is sharp up to the
constant, since there are groups with pinched growth $a^r$
where~(\ref{e.ii}) cannot be improved.  Examples include the
lamplighter on $\Z$~\cite{LPP:lamplighter}. 
See \cite{Ha} for a discussion of growth rates of groups and
many related open problems. 

\medskip

Regarding pinched polynomial growth, it is known that for every $d>1$
there is a tree with pinched growth $r^d$ containing arbitrarily large
sets $A$ with $|\p A|=1$, see, e.g., \cite{BS}.
\medskip

\begin{prob}\label{pr.tree}
Does every graph of a pinched exponential growth
contain a tree with pinched exponential growth? 
\end{prob}

In~\cite{BS:Cheeger} it was shown that every graph satisfying the linear isoperimetric
inequality $|\p A|\ge c\,|A|$ ($c>0$) contains a tree satisfying such an inequality,
possibly with a different constant.  The question whether one can find a {\bf spanning}
tree  with a linear isoperimetric inequality was asked earlier~\cite{DSS}.
It follows from Theorem~\ref{main} that
a tree with pinched exponential growth satisfies the linear isoperimetric
inequality.  (If a tree satisfies $|\p A| \ge 3$ for every vertex set $A$
of size at least $k$, then every path of $k$ vertices in the tree must contain
a branch point, a point whose removal will give at least $3$ infinite components.
Consequently the tree contains a modified infinite
binary tree, where every edge is subdivided into at most $k$ edges.)
Consequently, Problem~\ref{pr.tree} is equivalent to the question whether
every graph with pinched exponential growth contains a tree satisfying a linear
isoperimeteric inequality.

\bigskip

As a warm up for the proof of Theorem~\ref{main},
here is an easy argument showing that when $G$ has pinched growth $a^r$ it satisfies
a two-dimensional isoperimetric inequality.  Let
$A\subset V(G)$ be finite. Let $v$ be a vertex of A that is
farthest from $\partial A$, and let $r$ be the distance from $v$
to $\partial A$ . Note that $B(v,2r)\subset \bigcup_{u \in \partial
A} B(u,r)$. This gives,
   $a^{2r} \le O(1)\,   |\partial A |\,  a^r$, and therefore
  $O(1)\,|\partial A | \ge  a^r$. On the other hand,
$\bigcup_{u \in  \partial A } B(u,r)\supset A$, which gives
$O(1)\, |\partial A|\,   a^r \ge |A|$. Hence, $O(1)\,   |\partial A|^2 \ge |A|$.

\bigskip
{\noindent \bf Proof of Theorem.}
For vertices $v,u$ set $z(v,u) := a^{-d(v,u)}$, where $d(v,u)$ is the graph distance
between $v$ and $u$ in $G$.

We estimate in two ways the quantity
$$
Z=Z_A := \sum_{v \in A} \sum_{u \in \partial
A} z(v,u)\,.
$$

Fix $v \in A$.
For every $w \notin A$, fix some geodesic path from $v$ to $w$, and let $w'$ be the first
vertex in $\partial A$ on this path.  Let $R$ be sufficiently large so that
$\bl|B(v,R)\br|\ge 2\,|A|$, and set $W:=B(v,R)\setminus A$.  Then
$$
\bl| \{ (w,w') : w \in W \} \br|  =  | W| \ge  a^R/O(1)\,.
$$
On the other hand, we may estimate the left hand side by
considering all possible $u \in \partial A$ as candidates for
$w'$.  If $w\in W $, then $d(v,w')+d(w',w)\le R$.
Thus, each $u$ is equal to $w'$ for at most $O(1)\, a^{R-d(v,u)}$ vertices $w\in W$.
This gives
$$
      \bl| \{ (w,w') : w \in W \} \br|   \le  O(1) \sum_{u \in \partial A}
      a^{R-d(u,v)} = O(1)\,a^R\sum_{u\in\p A} z(v,u)\,.
$$
Combining these two estimates yields $O(1)\sum_{u\in \p A} z(v,u) \ge 1$.
By summing over $v$, this implies
\begin{equation}\label{e.Z1}
O(1)\,Z \ge   |A|\,.
\end{equation}

Now fix $u \in \partial A$,
set $m_r:=\bl|\{v\in A:d(v,u)=r\}\br|$, and consider
\begin{equation}\label{e.Zu}
Z(u):=\sum_{v \in A} z(v,u)=\sum_r m_r\,a^{-r}\,.
\end{equation}
For $r\le \log|A|/\log a$, we use the inequality  $m_r\le \bl|B(u,r)\br|=O(1)\,a^r$,
while for $r>\log|A|/\log a$, we use $m_r\le |A|$.
We apply these estimates to~(\ref{e.Zu}), and get
$Z(u)\le O(1)\,\log(2+|A|)$,
which gives $Z=\sum_{u\in\p A} Z(u) \le O(1)\, |\partial A|\, \log(2+ |A|)$.
Together with~(\ref{e.Z1}), this gives~(\ref{e.ii}).
\qed

\bigskip
Next, we present a slightly different
version of Theorem~\ref{main}, which also applies to finite graphs.

\begin{theo}\label{finite}
Let $G$ be a finite or infinite graph, $c>0$, $a>1$, $R\in\N$, and
suppose that $c^{-1}\,a^r\le |B(v,r)|\le c\,a^r$ holds
for all $r=1,2,\dots,R$ and for all $v\in V(G)$.
Then there is a constant $C=C(a,c)$, depending only on $a$ and $c$,
such that
$$
C\,\bl|\partial A\br| \ge |A| / \log(2+ |A|)
$$
holds for every finite $A\subset V(G)$ with $|A|\le C^{-1}\, a^R$.
\end{theo}

The proof is the same.  A careful inspection of the proof shows
that one only needs the inequality
$c^{-1}\,a^r\le |B(v,r)|$ to be valid for $v\in A$ and the
inequality
$|B(v,r)|\le c\,a^r$ only for $v\in\p A$.

\bigskip

\noindent {Acknowledgements:} We thank Thierry Coulhon and Iftach
Haitner  for useful disscusions.


\begin{thebibliography}{99}



\bibitem{BaS} L. Babai and M. Szegedy,
Local expansion of symmetrical graphs, Combin. Probab. Comput. 1
(1992), no. 1, 1--11.


\bibitem{BS:Cheeger}
I. Benjamini and O. Schramm,
Every graph with a positive Cheeger constant contains a tree with a positive Cheeger constant.
Geom. Funct. Anal. vol. 7 (1997), no. 3, 403--419.

\bibitem{BS}
I. Benjamini and O. Schramm, Recurrence of distributional limits
of finite planar graphs. Electron. J. Probab.  6  (2001), no. 23,
13 pp. (electronic).

\bibitem{CS}
T. Coulhon and L. Saloff-Coste, Isop\'erim\'etrie pour les groupes
et les vari\'et\'es. (French) [Isoperimetry for groups and
manifolds]  Rev. Mat. Iberoamericana  9  (1993),  no. 2, 293--314.

\bibitem{DSS}
W. Deuber, M. Simonovits and V. S\'os, A note on paradoxical
metric spaces. Studia Sci. Math. Hungar. 30 (1995), no. 1-2,
17--23.



\bibitem{Ha}
 P. de la Harpe,
 Topics in geometric group theory. Chicago Lectures in Mathematics.
 University of Chicago Press, Chicago, IL, 2000. vi+310 pp.

\bibitem{LPP:lamplighter}
R. Lyons, R. Pemantle and Y. Peres, Random walks on the lamplighter group. Ann. Probab. 24 (1996), no. 4, 1993--2006.


\end{thebibliography}
\end{document}